\documentclass[12pt,thmsa]{article}
\usepackage{amsfonts}
\usepackage{myart}

\normalbaselines
\font\tenbm=cmmib10
\font\sevenbm=cmmib7
\font\fivebm=cmmib5
\newfam\bmfam
\textfont\bmfam=\tenbm \scriptfont\bmfam=\sevenbm
\scriptscriptfont\bmfam=\fivebm
{\count0=\number\bmfam \multiply\count0 by "100
\def\defbgreek#1#2#3{{\count1=\count0 \advance\count1 by "#2#3
  \global\mathchardef#1=\count1 }}
\defbgreek\balpha  0B \defbgreek\brho       1A
\defbgreek\bbeta   0C \defbgreek\bsigma     1B
\defbgreek\bgamma  0D \defbgreek\btau       1C
\defbgreek\bdelta  0E \defbgreek\bupsilon   1D
\defbgreek\bepsilon0F \defbgreek\bphi       1E
\defbgreek\bzeta   10 \defbgreek\bchi       1F
\defbgreek\bmeta   11 \defbgreek\bpsi       20
\defbgreek\btheta  12 \defbgreek\bomega     21
\defbgreek\biota   13 \defbgreek\bvarepsilon22
\defbgreek\bkappa  14 \defbgreek\bvartheta  23
\defbgreek\blambda 15 \defbgreek\bvarpi     24
\defbgreek\bmu     16 \defbgreek\bvarrho    25
\defbgreek\bnu     17 \defbgreek\bvarsigma  26
\defbgreek\bxi     18 \defbgreek\bvarphi    27
\defbgreek\bpi     19}

\newtheorem{opred}{Definition}[section]
\newtheorem{remark}{Remark}[section]
\newtheorem{theorem}{Theorem}[section]

\begin{document}

\author{Yuri A. Rylov}
\title{Geometry without Topology}
\date{Institute for Problems in Mechanics, Russian Academy of Sciences, \\
101-1, Vernadskii Ave., Moscow, 117526, Russia. \\
email: rylov@ipmnet.ru\\
Web site: {$http://rsfq1.physics.sunysb.edu/\symbol{126}rylov/yrylov.htm$}.\\
or mirror Web site: {$http://194.190.131.172/\symbol{126}rylov/yrylov.htm$}.}
\maketitle

\begin{abstract}
The proper Euclidean geometry is considered to be metric space and described
in terms of only metric and finite metric subspaces ($\sigma $-immanent
description). Constructing the geometry, one does not use topology and
topological properties. For instance, the straight, passing through points $%
A $ and $B$, is defined as a set of such points $R$ that the area $S(A,B,R)$
of triangle $ABR$ vanishes. The triangle area is expressed via metric by
means of the Hero's formula, and the straight appears to be defined only via
metric, i.e. without a reference to (topological) concept of curve.
(Usually, the straight is defined as the shortest curve, connecting two
points $A$ and $B$). Such a construction of geometry is free from such
restrictions as continuity and dimensionality of the space which are
generated by a use of topology but not by the geometry in itself. At such a
description all information on the geometry properties (such as uniformity,
isotropy, continuity and degeneracy) is contained in metric. The Riemannian
geometry is constructed by two different ways: (1) by conventional way on
the basis of metric tensor, (2) as a result of modification of metric in the
$\sigma $-immanent description of the proper Euclidean geometry. The two
obtained geometries are compared. The convexity problem in geometry and the
problem of collinearity of vectors at distant points are considered. The
nonmetric definition of curve is shown to be a concept of the proper
Euclidean geometry which is inadequate to any non-Euclidean geometry
\end{abstract}

\section{Introduction}

There are several methods of the proper Euclidean geometry\footnote{%
We use the term ''Euclidean geometry'' as a collective concept with respect
to terms ''proper Euclidean geometry'' and ''pseudo-Euclidean geometry''. In
the first case the eigenvalues of the metric tensor matrix have similar
signs, in the second case they have different signs.} description. The most
old way of description is the axiomatic conception of the Euclidean
geometry. The proper Euclidean geometry is described in terms of points,
straights and planes, which are determined by their properties in terms of
axioms. Some axioms describe properties of natural geometric objects
(points, straights and planes), other axioms describe such properties of
proper Euclidean geometry as uniformity, isotropy, continuity and degeneracy%
\footnote{%
In general case the set of vectors, parallel to the given vector, forms a
cone. Degeneracy of the proper Euclidean geometry means that this cone
degenerates into a line.}.

Real geometry of the space-time is uniform only approximately, and one needs
to consider the geometries which should not be uniform, isotropic,
continuous and degenerate. In other words, one needs to generalize and
modify the proper Euclidean geometry. But it is quite impossible one to
modify axioms of the proper Euclidean geometry, and one needs to describe
the proper Euclidean geometry in the form, containing numerical
characteristics which may be modified rather easily. Such numerical
characteristic of the proper Euclidean geometry is the metric $\rho (P,Q)$,
describing distance between any two points $P$ and $Q$ of the space. After
modification of the Euclidean metric a new geometry appears, which may have
other properties than uniformity, isotropy, continuity and degeneracy.

Usually for construction of (Riemannian) geometry one uses the following
logical scheme

\begin{equation}
\begin{array}{ccccccc}
\begin{array}{c}
\mbox{coordinate} \\
\mbox{system}
\end{array}
& \longrightarrow &
\begin{array}{c}
\mbox{infinitesimal} \\
\mbox{distance}
\end{array}
& \longrightarrow &
\begin{array}{c}
\mbox{set of } \\
\mbox{geodesics}
\end{array}
&
\begin{array}{c}
\nearrow \\
\searrow
\end{array}
&
\begin{array}{c}
\mbox{geometry} \\
\\
\begin{array}{c}
\mbox{finite} \\
\mbox{distance}
\end{array}
\end{array}
\end{array}
\label{a.1}
\end{equation}

As it follows from this scheme for construction of geometry one needs a
coordinate system and a system of geodesics. The coordinate system is
necessary for introduction of infinitesimal distance. The geodesic is
defined as a shortest curve (line), connecting two points. Thus, the
considered construction of geometry refers to the concept of a curve. The
curve is a topological object, defined as a continuous mapping of the real
axis onto the geometrical space of points. As a result the topology is
considered usually to be a necessary element of geometry. According to (\ref
{a.1}) one cannot construct geometry without a use of topology (in the form
of a curve). Actually the topology is only a mathematical tool, using for
construction of geometry. To prove this, it is sufficient to construct
geometry without a reference to the topological concept of a curve. We shall
make this in the present paper.

The geometry is constructed in accord with the following logical scheme
\begin{equation}
\begin{array}{cc}
\begin{array}{c}
\mbox{finite} \\
\mbox{distance}
\end{array}
&
\begin{array}{cc}
\begin{array}{c}
\nearrow \\
\\
\searrow
\end{array}
&
\begin{array}{c}
\mbox{geometry} \\
\\
\begin{array}{c}
\mbox{set of } \\
\mbox{geodesics}
\end{array}
\end{array}
\end{array}
\end{array}
\label{a.2}
\end{equation}
where geometry is constructed independently of a possibility of the
geodesics construction. It is possible such a situation, when the geometry
can be constructed, whereas geodesics (the shortest curves) cannot. Such a
situation is not exotic, because the real space-time geometry appears to be
of such a kind. Timelike geodesics of the space-time are substituted by thin
hallow tubes. Thickness of the tubes is microscopic. Describing macroscopic
phenomena, one may neglect the tube thickness and substitute the tubes by
lines. Then geometry may be considered as a degenerate one (the tubes
degenerate into lines). Describing microscopic phenomena, one may not
neglect the thickness of tubes, because the thickness of tubes
(nondegeneracy of geometry) is a reason of quantum effects. Besides the
geometry constructed in accord with the scheme (\ref{a.2}) is free from such
constraints as continuity and degeneracy, imposed by a use of the concept of
a curve.

To carry out the idea of nondegenerate geometry, let us give some
definitions which help us to formulate the problem of generalization and
modification of the proper Euclidean geometry.

\begin{opred}
\label{d0} The metric space $M=\{\rho ,\Omega \}$ is a set $\Omega $ of
points $P\in \Omega $ with the metric $\rho $ given on $\Omega \times \Omega
$
\begin{equation}
\rho :\quad \Omega \times \Omega \rightarrow D_{+}\subset \Bbb{R}
\label{a1.1}
\end{equation}
\begin{equation}
\rho (P,P)=0,\qquad \rho (P,Q)=\rho (Q,P),\qquad \forall P,Q\in \Omega
\label{a1.2}
\end{equation}
\begin{equation}
D_{+}=[0,\infty ),\qquad \rho (P,Q)=0,\quad \mbox{if and only if}\mathrm{%
\;\;\;\;}P=Q,\qquad \forall P,Q\in \Omega  \label{a1.3}
\end{equation}
\begin{equation}
\rho (P,Q)+\rho (Q,R)\geq \rho (P,R),\qquad \forall P,Q,R\in \Omega
\label{a1.4}
\end{equation}
\end{opred}

\begin{opred}
\label{d1} Any subset $\Omega ^{\prime }\subset \Omega $ of points of the
metric space $M=\{\rho ,\Omega \}$, equipped with the metric $\rho ^{\prime
} $ which is a contraction $\rho |_{\Omega ^{\prime }\times \Omega ^{\prime
}}$ of the mapping (\ref{a1.1}). on the set $\Omega ^{\prime }\times \Omega
^{\prime }$ is called the metric subspace $M^{\prime }=\{\rho ^{\prime
},\Omega ^{\prime }\}$ of the metric space $M=\{\rho ,\Omega \}$.
\end{opred}

It is easy to see that the metric subspace $M^{\prime }=\{\rho ^{\prime
},\Omega ^{\prime }\}$ is a metric space.

\begin{opred}
\label{d2} The metric space $M=\{\rho ,\Omega \}$ is called finite, if the
set $\Omega $ contains a finite number of points. The finite metric subspace
$M(\mathcal{P}^n)=\{\rho ,\mathcal{P}^n\}$ of $M=\{\rho ,\Omega \}$,
consisting of $n+1$ points $\mathcal{P}^n\equiv \{P_0,P_1,\ldots
P_n\}\subset \Omega ,\quad $ $n=0,1,\ldots $is called the $n$th order metric
subspace.
\end{opred}

The proper Euclidean space may be considered to be a kind of metric space $%
E=\left\{ \rho _E,\Omega \right\} $. Being a metric space, the proper
Euclidean space and geometry on this space can be described in terms of only
metric $\rho $ and of finite metric subspaces. The finite metric subspaces $%
M(\mathcal{P}^n)$ are the simplest constituents of the metric space. Some
properties of finite metric subspaces $M(\mathcal{P}^n)$ were investigated
by Blumenthal \cite{B53}, but he did not consider them to be primitive
fundamental objects of metric space as we do. Metric space $M(\mathcal{P}^n)$%
, consisting of $n+1$ points and having nonvanishing length (concept of the
length will be defined further), generates in the proper Euclidean space $n$%
-dimensional plane $\mathcal{L}_n(\mathcal{P}^n)$, which appears to be an
attribute of $M(\mathcal{P}^n)$ and can be defined in terms of $M(\mathcal{P}%
^n)$.

\begin{opred}
\label{d1.6} Elementary geometrical object is a set of points having some
metric property.
\end{opred}

\begin{opred}
\label{d1.7} Geometrical object is a set of points derived as joins and
intersections of elementary geometrical objects.
\end{opred}

\noindent In other words, a geometrical object is a metric subspace $%
M_G=\{\rho,G\}$, $G\subset\Omega $ of metric space $M=\{\rho ,\Omega \}$.

\begin{opred}
\label{d1.8} Geometry is a totality of all propositions (definitions, axioms
and theorems) on properties of geometrical objects.
\end{opred}

\noindent In other words, the geometry is a totality of all propositions on
properties of all metric subspaces of the metric space $M=\{\rho ,\Omega \}$.

Let us consider some examples of elementary geometrical objects.

\begin{opred}
\label{d1.9} The sphere $\mathcal{S}\left( O;P\right) $, having its center
at the point $O$ and passing through the point $P$, is the set of points $%
R\in \Omega $ of the metric space $M=\{\rho ,\Omega \},$ defined by the
relation
\[
\mathcal{S}\left( O;P\right) =\left\{ R|\rho \left( O,R\right) =\rho \left(
O,P\right) \right\} ,\qquad O,P,R\in \Omega
\]
\end{opred}

The basic points $O$ and $P$, determining the sphere $\mathcal{S}\left(
O;R\right) $, are not equivalent, because $\mathcal{S}\left( O;R\right) $
and $\mathcal{S}\left( R;O\right) $ are different elementary geometrical
objects (different spheres). In particular, $P\in \mathcal{S}\left(
O;R\right) $, but $O\notin \mathcal{S}\left( O;R\right) $. The sphere $%
\mathcal{S}\left( O;P\right) $ is an attribute of zeroth order metric
subspaces $M(O)$ and $M(P)$ (or two points $O,P$).

\begin{opred}
\label{d1.10a} The circle cylinder $\mathcal{C}\left( P_1,P_2;P\right) $,
passing through the point $P$, with axis, determined by the basic points $%
P_1,P_2$, is the set of points $R\in \Omega $ of the metric space $M=\{\rho
,\Omega \}$, defined by the relation
\begin{eqnarray*}
\mathcal{C}\left( P_1,P_2;P\right) =\left\{ R|S_2\left( P_1,P_2,R\right)
=S_2\left( P_1,P_2,P\right) \right\} , \\
P_1,P_2,P,R\in \Omega
\end{eqnarray*}
\end{opred}

\noindent where $S_2\left( P_1,P_2,R\right) $ is the area of the triangle
with vertices at the points $P_1,P_2,R$. If the areas of triangles $%
\triangle P_1P_2R$ and $\triangle P_1P_2P$ are equal, the heights (radii)
dropped from the vertices $R$ and $P$ of these triangles onto their common
base $P_1P_2$ (axis of the cylinder) are also equal. The triangle area $%
S_2\left( P_1,P_2,R\right) $ can be expressed via metric by the Hero's
formula
\[
S_2\left( A,B,C\right) =\sqrt{p\left( p-a\right) \left( p-b\right) \left(
p-c\right) },
\]
\[
a=\rho \left( B,C\right) ,\qquad b=\rho \left( A,C\right) ,\qquad c=\rho
\left( A,B\right) ,\qquad p=\left( a+b+c\right) /2
\]
The circle cylinder $\mathcal{C}\left( P_1,P_2;P\right) $ is an attribute of
two finite metric subspaces $M(P_1,P_2)$ and $M(P)$.

\begin{opred}
\label{d1.10} The ellipsoid $\mathcal{E}\left( P_1,P_2;P\right) $, having
its focuses at the basic points $P_1,P_2$ and passing through the point $P$
is the set of points $R\in \Omega $ of the metric space $M=\{\rho ,\Omega \}$%
, defined by the relation
\begin{eqnarray*}
\mathcal{E}\left( P_1,P_2;P\right) =\left\{ R|\rho \left( P_1,R\right) +\rho
\left( P_2,R\right) =\rho \left( P_1,P\right) +\rho \left( P_2,P\right)
\right\} , \\
P_1,P_2,P,R\in \Omega
\end{eqnarray*}
\end{opred}

The ellipsoid $\mathcal{E}\left( P_1,P_2;P\right) $ is an attribute of two
finite metric subspaces $M(P_1,P_2)$ and $M(P)$. If $P_1\neq P$ and $P_2\neq
P$, the points $P_1,P_2\notin \mathcal{E}\left( P_1,P_2;P\right) ,$ but the
point $P\in \mathcal{E}\left( P_1,P_2;P\right) $. If the point $P=P_1$, the
ellipsoid $\mathcal{E}\left( P_1,P_2;P\right) $ degenerates into segment $%
\mathcal{T}_{[P_1P_2]}$ between the points $P_1$ and $P_2$ of the straight
line $\mathcal{T}_{P_1P_2}$, passing through the points $P_1$ and $P_2$. The
segment $\mathcal{T}_{[P_1P_2]}$ is defined as follows.

\begin{opred}
\label{d1.11} The segment $\mathcal{T}_{[P_1P_2]}$ of the straight between
the basic points $P_1,P_2$ is the set of points $R\in \Omega $ of the metric
space $M=\{\rho ,\Omega \},$ defined by the relation
\begin{equation}
\mathcal{T}_{[P_1P_2]}=\left\{ R|\rho \left( P_1,R\right) +\rho \left(
P_2,R\right) -\rho \left( P_1,P_2\right) =0\right\} ,\qquad P_1,P_2,R\in
\Omega  \label{b1.1}
\end{equation}
\end{opred}

The segment $\mathcal{T}_{[P_1P_2]}$ is an elementary geometrical object
which does not depend on the order of points $P_1,P_2$. Besides both basic
points $P_1,P_2\in \mathcal{T}_{[P_1P_2]}.$ The segment $\mathcal{T}%
_{[P_1P_2]}$ is an attribute of the first order metric subspace $M\left(
P_1,P_2\right) $ in the sense that $\mathcal{T}_{[P_1P_2]}$ is determined by
the metric subspace $M\left( P_1,P_2\right) $ in itself. For instance, the
sphere $\mathcal{S}\left( O;P\right) $ is determined by the points $O,P$ of
the metric subspace $M\left( O,P\right) ,$ but not by the metric subspace $%
M\left( O,P\right) $ in itself, and the sphere $\mathcal{S}\left( O;P\right)
$ is not an attribute of the metric subspace $M\left( O,P\right) $, but it
is an attribute of two zeroth order metric subspaces $M(O)$ and $M(P)$ (or
two points $O,P$).

\begin{opred}
\label{d1.12} The elemetary geometrical object which is an attribute of the $%
n$th order metric subspace $M\left( \mathcal{P}^n\right) $ is the $n$th
order natural geometric object (the $n$th order NGO).
\end{opred}

Such geometrical objects as a point, an Euclidean straight, and an Euclidean
plane are NGOs of the proper Euclidean geometry. The point $P_0$ is the
zeroth order NGO $\mathcal{T}_{P_0}$ of the proper Euclidean geometry which
is determined by the zeroth order metric subspace $M\left( P_0\right) =P_0$.
The straight $\mathcal{T}_{P_0P_1}$ of the proper Euclidean geometry is the
first order NGO which is determined by the first order metric subspace $%
M\left( P_0,P_1\right) $. It means, in particular, that $\mathcal{T}%
_{P_0P_1}=\mathcal{T}_{P_1P_0}$.

The two-dimensional plane $\mathcal{T}_{P_0P_1P_2}$ of the proper Euclidean
geometry is the second order NGO, determined by the second order metric
subspace $M\left( P_0,P_1,P_2\right) $. It means that the NGO $\mathcal{T}%
_{P_0P_1P_2}$ does not depend on the order of basic points $P_0,P_1,P_2$,
which determine $\mathcal{T}_{P_0P_1P_2}$. It does not always happen that
the second order metric subspace $M\left( P_0,P_1,P_2\right) $ determines $%
\mathcal{T}_{P_0P_1P_2}$. Only $M\left( P_0,P_1,P_2\right) \not \subset
\mathcal{T}_{P_0P_1}$ enables to determine $\mathcal{T}_{P_0P_1P_2}$.

For explicit determination of the $n$th order NGO one needs to attribute a
length $|M\left( \mathcal{P}^n\right) |$ to any $n$th order metric subspace $%
M\left( \mathcal{P}^n\right) $

\begin{opred}
\label{d1.13} The squared length $\left| M\left( \mathcal{P}^n\right)
\right| ^2 $of the $n$th order metric subspace $M\left( \mathcal{P}^n\right)
\subset \Omega $ of the proper Euclidean space $E=\left\{ \rho _E,\Omega
\right\} $ is the real number.
\[
\left| M\left( \mathcal{P}^n\right) \right| ^2=\left( n!S_n(\mathcal{P}%
^n)\right) ^2=F_n\left( \mathcal{P}^n\right)
\]
where $S_n(\mathcal{P}^n)$ is the Euclidean volume of the $(n+1)$-edr with
vertices at points $\mathcal{P}^n\equiv \{P_0,P_1,\ldots P_n\}\subset \Omega
$.
\end{opred}

In the proper Euclidean geometry the volume $S_n(\mathcal{P}^n)$ of the $%
(n+1)$-edr and the value $F_n\left( \mathcal{P}^n\right) $ of the function $%
F_n$, connected with it, can be expressed in terms of metric $\rho $ by
means of relations
\begin{equation}
F_n:\quad \Omega ^{n+1}\rightarrow \Bbb{R},\qquad \Omega
^{n+1}=\bigotimes\limits_{k=1}^{n+1}\Omega ,\qquad n=1,2,\ldots  \label{a1.5}
\end{equation}
\begin{equation}
F_n\left( \mathcal{P}^n\right) =\det ||\left( \mathbf{P}_0\mathbf{P}_i.%
\mathbf{P}_0\mathbf{P}_k\right) ||,\qquad P_0,P_i,P_k\in \Omega ,\qquad
i,k=1,2,...n  \label{a1.6}
\end{equation}
\begin{eqnarray}
\left( \mathbf{P}_0\mathbf{P}_i.\mathbf{P}_0\mathbf{P}_k\right) &\equiv
&\Gamma \left( P_0,P_i,P_k\right) \equiv \sigma \left( P_0,P_i\right)
+\sigma \left( P_0,P_k\right) -\sigma \left( P_i,P_k\right) ,  \label{a1.7}
\\
i,k &=&1,2,...n.  \nonumber
\end{eqnarray}
where the function $\sigma $ is defined via metric $\rho $ by the relation
\begin{equation}
\sigma (P,Q)\equiv \frac 12\rho ^2(P,Q),\qquad \forall P,Q\in \Omega .
\label{a1.8}
\end{equation}
and $\mathcal{P}^n$ denotes $n+1$ points $P_0,P_1,\ldots ,P_n$ of $\Omega $%
\begin{equation}
\mathcal{P}^n=\{P_0,P_1,\ldots ,P_n\}\subset \Omega  \label{a1.9}
\end{equation}
The function $\sigma ,$ called world function \cite{S60}, is very important
quantity which may be used instead of metric $\rho $. In many cases a use of
the function $\sigma $ appears to be more convenient than a usage of metric $%
\rho $. The squared length $\left| M\left( \mathcal{P}^n\right) \right|
^2=F_n\left( \mathcal{P}^n\right) $ is calculated for the proper Euclidean
space, but the expression (\ref{a1.6}) - (\ref{a1.8}) may be used for any
finite subspaces of any metric space, because it contains only world
function $\sigma $ (metric $\rho )$ and may be calculated for any metric
space.

\begin{opred}
\label{d3} A description is called $\sigma $-immanent, if it does not
contain any references to objects or concepts other than finite subspaces of
the metric space and its metric.
\end{opred}

Prefix $\sigma $ in the term ''$\sigma $-immanent'' associates with the
world function $\sigma $. Concept of $\sigma $-immanent description is very
important for modification of the proper Euclidean geometry. Considering the
proper Euclidean geometry to be a standard geometry and defining a
geometrical object there in a $\sigma $-immanent way, one can use this
definition in any metric space.

Note that definition of geometrical objects is a principal problem of the
metric geometry, i.e. the geometry, generated by the metric space. The
shortest (line), connecting two arbitrary points $P,Q\in \Omega $ of the
metric space $\{\rho ,\Omega \}$, is the basic geometrical object which is
constructed usually in the metric space \cite{ABN86}. One can construct an
angle, triangle, different polygons from segments of the shortest.
Construction of two-dimensional and three-dimensional planes in the metric
space is rather problematic. At any rate it is unclear how one could
construct these planes, using the shortest as the main geometrical object. A
possibility of the metric space description in terms of only the shortest is
restricted. Although exhibiting ingenuity, such a description may be
constructed. For instance, A.D.~Alexandrov showed that internal geometry of
two-dimensional boundaries of convex three-dimensional bodies may be
represented in terms of metric \cite{A48}. Apparently, without introducing
geometric objects which are analogs of two-dimensional plane, the solution
of similar problem for three-dimensional boundaries of four-dimensional
bodies is very difficult.

Note that constraints (\ref{a1.3}), (\ref{a1.4}), imposed on metric, are
necessary only for constructing the shortest. The shortest, determined by
two points $P_1,P_2$, may be replaced by the $\sigma $-immanent definition (%
\ref{b1.1}) of segment $\mathcal{T}_{[P_1P_2]}$, which coincides with the
shortest in the metric space, described by the definition \ref{d0}. This
definition in itself does not need constraint (\ref{a1.3}), describing
definiteness of the metric space, and constraint (\ref{a1.4}), describing
one-dimensionality of the segment $\mathcal{T}_{[P_1P_2]}.$ If the metric is
not restricted by constraint (\ref{a1.4}), the segment $\mathcal{T}%
_{[P_1P_2]}$ takes the shape of a hallow tube, reminding ellipsoid,
described by definition \ref{d1.10}. If the constraint (\ref{a1.4}) is
strengthened ($\leq $ is replaced by $<$), the segment $\mathcal{T}%
_{[P_1P_2]}$ degenerates into two points $P_1,P_2$. The case of the
one-dimensional shortest is intermediate between the two cases.

In the case of the proper Euclidean space, considered to be a metric space,
the first order NGO, defined by (\ref{b1.1}) is one-dimensional line. It is
not clear whether one-dimensionality is a special property of the Euclidean
geometry, or it is a property of any geometry in itself. We do not see, why
one should insist on the one-dimensionality of the first order NGO $\mathcal{%
T}_{[P_1P_2]}$ in the case of an arbitrary modification of the proper
Euclidean geometry. First, it is useful to consider the most general
modification of the proper Euclidean geometry. Second, at the end of
investigation, if it appears to be necessary, one can always reduce a degree
of generalization, imposing additional constraints.

In the proper Euclidean space the $n$-dimensional plane ($n$th order NGO) $%
n=1,2,\ldots $ is defined as follows

\begin{opred}
\label{d1.14} \quad The $n$th order metric subspace $\;M\left( \mathcal{P}%
^n\right) \;$ of unvanishing $\;$ length $\;$ $\left| M\left( \mathcal{P}%
^n\right) \right| ^2=F_n\left( \mathcal{P}^n\right) \neq 0$ determines the $%
n $th order tube (the $n$th order NGO) $\mathcal{T}\left( \mathcal{P}%
^n\right) $ by means of the relation
\begin{equation}
\mathcal{T}\left( \mathcal{P}^n\right) \equiv \mathcal{T}_{\mathcal{P}%
^n}=\left\{ P_{n+1}|F_{n+1}\left( \mathcal{P}^{n+1}\right) =0\right\}
,\qquad P_i\in \Omega ,\qquad i=0,1\ldots n+1,  \label{b1.3}
\end{equation}
where the function $F_n$ is defined by the relations (\ref{a1.5}), (\ref
{a1.7})
\end{opred}

The $n$th order tube $\mathcal{T}_{\mathcal{P}^n}$ which is an analog of the
$n$-dimensional Euclidean plane may be constructed in any metric space, as
far as its definition \ref{d1.14} is $\sigma $-immanent. It may be defined
also in the metric space with omitted constraints (\ref{a1.3}), (\ref{a1.4}%
), imposed usually on the metric. We shall refer to such a generalized
metric space as the $\sigma $-space. The geometry, generated by the $\sigma $%
-space, will be referred to as T-geometry (tubular geometry).

\begin{opred}
\label{d3.1.1} $\sigma $-space $V=\{\sigma ,\Omega \}$ is nonempty set $%
\Omega $ of points $P$ with given on $\Omega \times \Omega $ real function $%
\sigma $
\begin{equation}
\sigma :\quad \Omega \times \Omega \to \Bbb{R},\qquad \sigma (P,P)=0,\qquad
\sigma (P,Q)=\sigma (Q,P)\qquad \forall P,Q\in \Omega .  \label{a2.1}
\end{equation}
\end{opred}

The function $\sigma $ is called world function, or $\sigma $-function. The
metric $\rho $ may be introduced in the $\sigma $-space by means of the
relation (\ref{a1.8}). If $\sigma $ is positive, metric $\rho $ is also
positive, but if $\sigma $ is negative, the metric is imaginary.

\begin{opred}
\label{d3.1.1ae}. Nonempty subset $\Omega ^{\prime }\subset \Omega $ of
points of the $\sigma $-space $V=\{\sigma ,\Omega \}$ with the world
function $\sigma ^{\prime }=\sigma |_{\Omega ^{\prime }\times \Omega
^{\prime }}$, which is a contraction of $\sigma $ on $\Omega ^{\prime
}\times \Omega ^{\prime }$, is called $\sigma $-subspace $V^{\prime
}=\{\sigma ^{\prime },\Omega ^{\prime }\}$ of $\sigma $-space $V=\{\sigma
,\Omega \}$.
\end{opred}

Further the world function $\sigma ^\prime = \sigma |_{\Omega ^{\prime
}\times\Omega ^\prime }$, which is a contraction of $\sigma $ will be
designed by means of $\sigma $. Any $\sigma$-subspace of $\sigma$-space is a
$\sigma$-space.

\begin{opred}
\label{d3.1.1ba}. $\sigma $-space $V=\{\sigma ,\Omega \}$ is called
isometrically embeddable in $\sigma $-space $V^{\prime }=\{\sigma ^{\prime
},\Omega ^{\prime }\}$, if there exists such a monomorphism $f:\Omega
\rightarrow \Omega ^{\prime }$, that $\sigma (P,Q)=\sigma ^{\prime
}(f(P),f(Q))$,\quad $\forall P,\forall Q\in \Omega ,\quad f(P),f(Q)\in
\Omega ^{\prime }$,
\end{opred}

Any $\sigma $-subspace $V^{\prime }$ of $\sigma $-space $V=\{\sigma ,\Omega
\}$ is isometrically embeddable in it.

\begin{opred}
\label{d3.1.1b}. Two $\sigma $-spaces $V=\{\sigma ,\Omega \}$ and $V^{\prime
}=\{\sigma ^{\prime },\Omega ^{\prime }\}$ are called to be isometric
(equivalent), if $V$ is isometrically embeddable in $V^{\prime }$, and $%
V^{\prime }$ is isometrically embeddable in $V$.
\end{opred}

\begin{opred}
\label{d2.2b} The $\sigma $-space $M=\{\rho ,\Omega \}$ is called a finite $%
\sigma $-space, if the set $\Omega $ contains a finite number of points.
\end{opred}

\begin{opred}
\label{d3.1.1bc}. The $\sigma $-subspace $M_n(\mathcal{P}^n)=\{\sigma ,%
\mathcal{P}^n\}$of the $\sigma $-space $V=\{\sigma ,\Omega \}$, consisting
of $n+1$ points $\mathcal{P}^n=\left\{ P_0,P_1,...,P_n\right\} $ is called
the $n$th order $\sigma $-subspace .
\end{opred}

All geometrical objects of T-geometry are obtained as follows. Geometrical
objects of the proper Euclidean geometry are defined in the $\sigma $%
-immanent form. Then they may be considered to be definitions of
corresponding geometrical objects in T-geometry. The world function $\sigma $
of the proper Euclidean space satisfies some $\sigma $-immanent relations,
describing special properties of the proper Euclidean geometry. Metric side
of these relations had been formulated and proved by Menger \cite{M28}.
Using our designations, we present this result in the form of theorem.

\begin{theorem}
\label{t1} The $\sigma $-space $V=\{\sigma ,\Omega \}$ is isomerically
embeddable in $n$-dimensional Euclidean space $E_n$, if and only if any $%
(n+2)$th order $\sigma $-subspace $M(\mathcal{P}^{n+2})\subset \Omega $ is
isometrically embeddable in $E_n$.
\end{theorem}

Unfortunately, the formulation of this theorem is not $\sigma $-immanent, as
far as it contains a reference to $n$-dimensional Euclidean space $E_n$
which is not defined $\sigma $-immanently. A more constructive version of
the $\sigma $-space Euclideness conditions is formulated in the form of the
following theorem.

\begin{theorem}
\label{c2}The $\sigma $-space $V=\{\sigma ,\Omega \}$ is the $n$-dimensional
Euclidean space, if and only if the following three $\sigma $-immanent
conditions are fulfilled.
\end{theorem}

\noindent I.
\begin{equation}
\exists \mathcal{P}^n\subset \Omega ,\qquad F_n(\mathcal{P}^n)\ne 0,\qquad
F_{n+1}(\Omega ^{n+2})=0,  \label{a3.4}
\end{equation}
II.
\[
\sigma (P,Q)={\frac 12}\sum_{i,k=1}^ng^{ik}(\mathcal{P}^n)[\Gamma
(P_0,P_i,P)-\Gamma (P_0,P_i,Q)]
\]
\begin{equation}
\times [\Gamma (P_0,P_k,P)-\Gamma (P_0,P_k,Q)],\qquad \forall P,Q\in \Omega
\label{a3.5}
\end{equation}
where $\Gamma (P_0,P_k,P)$ are defined by the relations (\ref{a1.7}). The
quantities $g^{ik}(\mathcal{P}^n),$ $(i,k=1,2,\ldots n)$ are defined by the
relations
\begin{equation}
\sum_{k=1}^ng_{ik}(\mathcal{P}^n)g^{kl}(\mathcal{P}^n)=\delta _i^l,\qquad
i,l=1,2,\ldots n  \label{a3.11}
\end{equation}
where
\begin{equation}
g_{ik}(\mathcal{P}^n)=\Gamma (P_0,P_i,P_k),\qquad i,k=1,2,\ldots n
\label{a3.9}
\end{equation}
III.\quad The relations
\begin{equation}
\Gamma (P_0,P_i,P)=x_i,\qquad x_i\in \Bbb{R},\qquad i=1,2,\ldots n,
\label{a3.12}
\end{equation}
considered to be equations for determination of $P\in \Omega $, have always
one and only one solution.

\begin{remark}
\label{r2} For the Euclidean space to be the proper Euclidean the
eigenvalues of the matrix $g_{ik}(\mathcal{P}^n)=\Gamma (P_0,P_i,P_k),\qquad
i,k=1,2,\ldots n$ are to be of the same sign, otherwise the Euclidean space
is pseudo-Euclidean.
\end{remark}

\begin{remark}
\label{r3} The condition (\ref{a3.4}) is a corollary of condition (\ref{a3.5}%
). It is formulated as a separate condition in order to separate definition
of dimension and that of the coordinate system.
\end{remark}

Let us note that all three conditions are written in $\sigma $-immanent
form. Proof of this theorem can be found in \cite{R00}. Now we consider how
results of this theorem can be used for construction of conventional
description of the proper Euclidean space in some rectilinear coordinate
system, starting from an abstract $\sigma $-space, satisfying conditions I -
III of the theorem.

Let there be $\sigma $-space $V=\{\sigma ,\Omega \},$ and it is known that
conditions I - III of the theorem are fulfilled. Then the $\sigma $-space $V$
is an Euclidean space, but the dimension $n$ of the space is unknown. To
determine the dimension $n$, let us take two different points $P_0,P_1\in
\Omega ,\;\;F_1(\mathcal{P}^1)=2\sigma (P_0,P_1)\neq 0$.

1. Let us construct the first order tube $\mathcal{T}\left( \mathcal{P}%
^1\right) $. If $\mathcal{T}\left( \mathcal{P}^1\right) =\Omega $, then
dimension of the $\sigma $-space $V$ \ $n=1$. If $\Omega \backslash \mathcal{%
T}\left( \mathcal{P}^1\right) \neq \emptyset ,\;\;\exists P_2\in \Omega
,\;\;P_2\notin \mathcal{T}\left( \mathcal{P}^1\right) ,$ and hence, $F_2(%
\mathcal{P}^2)\neq 0$.

2. Let us construct the second order tube $\mathcal{T}\left( \mathcal{P}%
^2\right) $. If $\mathcal{T}\left( \mathcal{P}^2\right) =\Omega $, then \ $%
n=2$, otherwise$\;\;\exists P_3\in \Omega ,\;\;P_3\notin \mathcal{T}\left(
\mathcal{P}^2\right) ,$ and hence, $F_3(\mathcal{P}^3)\neq 0$.

3. Let us construct the third order tube $\mathcal{T}\left( \mathcal{P}%
^3\right) $. If $\mathcal{T}\left( \mathcal{P}^3\right) =\Omega $, then $n=3$%
, otherwise$\;\;\exists P_4\in \Omega ,\;\;P_4\notin \mathcal{T}\left(
\mathcal{P}^3\right) ,$ and hence, $F_4(\mathcal{P}^4)\neq 0$.

4. Etc.

Continuing this process, one determines such $n+1$ points $\mathcal{P}^n$,
that the condition $\mathcal{T}\left( \mathcal{P}^n\right) =\Omega $ and,
hence, conditions (\ref{a3.4}) are fulfilled.

Then by means of relations
\begin{equation}
x_i\left( P\right) =\Gamma (P_0,P_i,P),\qquad i=1,2,\ldots n,  \label{a3.13}
\end{equation}
one attributes covariant coordinates $x\left( P\right) =\left\{
x_i(P)\right\} ,\;\;i=1,2,\ldots n$ to $\forall P\in \Omega $. Let $%
x=x\left( P\right) \in \Bbb{R}^n$ and $x^{\prime }=x\left( P^{\prime
}\right) \in \Bbb{R}^n.$ Substituting $\Gamma (P_0,P_i,P)=x$ and $\Gamma
(P_0,P_i,P^{\prime })=x_i^{\prime }$ in (\ref{a3.5}), one obtains the
conventional expression for the world function of the Euclidean space in the
rectilinear coordinate system
\begin{equation}
\sigma (P,P^{\prime })=\sigma _E(x,x^{\prime })={\frac 12}%
\sum_{i,k=1}^ng^{ik}(\mathcal{P}^n)\left( x_i-x_i^{\prime }\right) \left(
x_k-x_k^{\prime }\right)  \label{a3.14}
\end{equation}
where $g^{ik}(\mathcal{P}^n)$, defined by relations (\ref{a3.9}) and (\ref
{a3.11}), is the contravariant metric tensor in this coordinate system.

Condition III of the theorem states that the mapping
\[
x:\;\Omega \rightarrow \Bbb{R}^n
\]
described by the relation (\ref{a3.13}) is a bijection, i.e. for $\forall
y\in \Bbb{R}^n$ there exists such one and only one point $Q\in \Omega ,$
that $y=x\left( Q\right) $.

Thus, on the base of the world function, given on abstract set $\Omega
\times \Omega $, one can determine the dimension $n$ of the Euclidean space,
construct rectilinear coordinate system with the metric tensor $g_{ik}(%
\mathcal{P}^n)=\Gamma (P_0,P_i,P_k),\qquad i,k=1,2,\ldots n$ and describe
all geometrical objects which are determined in terms of coordinates. The
Euclidean space and Euclidean geometry is described in terms and only in
terms of world function (metric). Changing the world function, one obtains
another $\sigma $-space and another (non-Euclidean) geometry. One should
expect that another geometry is also described completely in terms of the
world function. The properties of geometrical objects may appear other than
the properties of these objects in the proper Euclidean geometry. For
instance, in the Euclidean geometry $\mathcal{T}_{P_0P_1}\subset \mathcal{T}%
_{P_0P_1P_2},$ i.e. the straight, passing through the points $P_0$ and $P_1$%
, belongs to any two-dimensional plane, passing through these points. To
prove these statement, one needs to use the relations (\ref{a3.5}). In the
case of non-Euclidean geometry the relation $\mathcal{T}_{P_0P_1}\subset
\mathcal{T}_{P_0P_1P_2}$ is invalid, in general.

Another example. Two circle cylinders $\mathcal{C}\left( P_0,P_1;P\right) $
and $\mathcal{C}\left( P_0,P_1^{\prime };P\right) ,$ $P_1^{\prime }\in
\mathcal{T}_{[P_0P_1]},\;P_1^{\prime }\neq P_1,$ $P_1^{\prime }\neq P_0$
coincide in the proper Euclidean geometry, but they are different
geometrical objects in non-Euclidean geometry.

In the proper Euclidean geometry there exists geometrical object called line.

\begin{opred}
\label{d.line}The broken line $\mathcal{T}_{\mathrm{br}}$ is the set of
connected straight segments $\mathcal{T}_{[P_iP_{i+1}]}$%
\begin{equation}
\mathcal{T}_{\mathrm{br}}=\bigcup\limits_i\mathcal{T}_{[P_iP_{i+1}]}
\label{a1.20}
\end{equation}
\end{opred}

The continuous line (or curve) is defined as a limit of the broken line $%
\mathcal{T}_{\mathrm{br}}$ at $P_i\rightarrow P_{i+1},$ $(i=0,\pm 1,\pm
2,\ldots )$. The smooth line is defined as a limit of (\ref{a1.20}) at $%
P_i\rightarrow P_{i+1},$ $(i=0,\pm 1,\pm 2,\ldots )$ under the constraint
that $\cos \angle P_{i-1}P_iP_{i+1}\rightarrow -1$. Defining the segment $%
\mathcal{T}_{[P_iP_{i+1}]}$ by means of definition \ref{d1.11}, one obtains
metric definition of broken line (\ref{a1.20}). To obtain metric definition
of the continuous line and that of smooth line, one needs to go to
corresponding limits in (\ref{a1.20}). According to this definition the line
is many-point geometrical object. This object is very complicated, because
their points are given independently (i.e. there are many degrees of
freedom).

On the other hand, in the proper Euclidean geometry there exists another
(nonmetric) definition of continuous line. The continuous line $\mathcal{L}$
is defined as a continuous mapping
\begin{equation}
L:\;I\rightarrow \Omega ,\qquad I=[0,1]\subset \Bbb{R}.  \label{b5.5}
\end{equation}
Strictly, the geometrical object is a set $\mathcal{L}=L(I)\subset\Omega$ of
points of the $\sigma$-space $V=\{\sigma ,\Omega\}$, but not the mapping (%
\ref{b5.5}) in itself. However, as far as the number set $I$ is fixed and
the same in all cases, then with some stipulations one can consider the
correspondence between the mapping $L$ and the set of images $\mathcal{L}%
=L(I)$ to be one-to-one. Then one can label the geometrical objects
(considered as $\sigma$-subspaces) by means of mappings (\ref{b5.5}) and
identify the mapping (\ref{b5.5}) with the geometrical object $\mathcal{L}$,
called curve (line).

In the proper Euclidean geometry the definition of line (\ref{b5.5}) agrees
with the definition \ref{a1.20}. But in non-Euclidean geometry definitions (%
\ref{b5.5}) and (\ref{a1.20}) do not agree, in general. Already in the
Riemannian geometry an application of definition (\ref{b5.5}) as one of
basic definitions poses problems.

In the Riemannian space the world function $\sigma _R(x,x^{\prime })$
between the points $x$ and $x^{\prime }$ is determined by the relation \cite
{S60}
\begin{equation}
\sigma _R(x,x^{\prime })=\frac 12\left( \int\limits_{_{\mathcal{L}%
_{[xx^{\prime }]}}}\sqrt{g_{ik}dx^idx^k}\right) ^2  \label{a3.15}
\end{equation}
where $\mathcal{L}_{[xx^{\prime }]}$ denotes segment of geodesic connecting
points $x$ and $x^{\prime }$. Let us use the world function (\ref{a3.15})
instead of the Euclidean world function (\ref{a3.14}) in the $\sigma $%
-immanent description of geometry. In other words, let us use for
construction of geometry the logical scheme (\ref{a.2}), but not (\ref{a.1}%
). One obtains the $\sigma $-Riemannian geometry which is expected to be
equivalent to the Riemannian geometry, because both the Riemannian geometry
and the $\sigma $-Riemannian one are two generalizations of the Euclidean
geometry, using the same world function which has to describe any geometry
completely. In reality, using for geometry construction different logical
schemes, the $\sigma $-Riemannian geometry and the Riemannian one coincide,
but not at all points,.

The point is that the world function is a fundamental object of the $\sigma $%
-Riemannian geometry, whereas it is a derivative object in the Riemannian
geometry, where the infinitesimal distance and the curve (line) are
fundamental objects. The line $\mathcal{L}$, defined by nonmetric definition
(\ref{b5.5}), is a complicated and fundamental structure of Riemannian
geometry, which is absent in such a form in the $\sigma $-Riemannian
geometry. The continuous line $\mathcal{L}$ in the $\sigma $-Riemannian
geometry may be defined as a limit of the broken tube (\ref{a1.20}). But it
is a derivative (not fundamental) geometrical object.

As a whole the situation looks as follows. The $\sigma $-Riemannian geometry
is constructed $\sigma $-immanently, i.e. on the base of metric and does not
need the nonmetric definition of line (\ref{b5.5}). The Riemannian geometry
is constructed on the base of infinitesimal metric $dS=\sqrt{g_{ik}dx^idx^k}$
(which coincide with the infinitesimal metric of the $\sigma $-Riemannian
geometry) and uses the nonmetric definition of line (\ref{b5.5}) for
definition of finite metric. As a result the finite metric of both
geometries coincide, but only in the whole domain $D=\Omega ,$ where both
geometries are defined. If one considers $\sigma $-Riemannian and Riemannian
geometries in some subdomain $D^{\prime }\subset D$, the finite metrics are
defined in $D^{\prime }$ in different ways for these geometries. For $\sigma
$-Riemannian geometry the finite metric in $D^{\prime }$ is defined as a
cotraction of the finite metric in $D$, whereas for Riemannian geometry the
finite metric is defined on the basis of system of geodesics inside $%
D^{\prime }$ which does not coincide, in general, with the system of
geodesics in $D$. The geodesic segment $\mathcal{L}_{[xx^{\prime }]}$ which
determines $\sigma _R(x,x^{\prime })$ is a lengthy geometrical object,
depending on the shape of the region $D^{\prime }$, where the Riemannian
geometry is defined. As a result the finite metrics of both geometries may
be different in $D^{\prime }\subset D$, although they coincide in $D$.

Note that the nonmetric definition of line (\ref{b5.5}) needs additional
constraints to be rather definite. Let us discuss these problems.

\section{Riemannian space and convexity problem}

The Riemannian space and the Riemannian geometry are introduced as follows. $%
n$-dimensional Riemannian space can be derived as a result of a
generalization of the $n$-dimensional proper Euclidean space, written in a
covariant form. Indeed, the $n$-dimensional Euclidean space $E_n=\left\{
\mathbf{g}_E,K,\Bbb{R}^n\right\} $ is described by the infinitesimal
distance written in the rectilinear coordinate system $K$
\begin{equation}
dS^2=g_{ik}dx^idx^k,\qquad g_{ik}\mathrm{=diag}\left\{ 1,1,...1\right\}
\label{p.1}
\end{equation}
$\mathbf{g}_E$ denotes the matrix $g_{ik}\mathrm{=diag}\left\{
1,1,...1\right\} $ of the metric tensor. In the arbitrary curvilinear
coordinate system $\tilde K$ the same distance have the form
\begin{equation}
dS^2=\tilde g_{ik}\left( \tilde x\right) d\tilde x^id\tilde x^k,\qquad \det {%
\tilde g_{ik}}\ne 0,  \label{p.2}
\end{equation}
Here $\tilde g_{ik}\left( x\right) $ is constrained by the relation
\begin{equation}
\tilde g_{ik}\left( x\right) =\sum\limits_{l=1}^n\frac{\partial f_l\left(
x\right) }{\partial x^i}\frac{\partial f_l\left( x\right) }{\partial x^k},
\label{p.3}
\end{equation}
where $f_l:\Bbb{R}^n\rightarrow \Bbb{R}$,\quad $l=1,2,\ldots n$ are $n$
functions restricted by one condition $\det ||\partial f_i/\partial
x^k||\neq 0,\quad i,k=1,2,...n$. If $\tilde g_{ik}\left( x\right) $ does not
satisfy the relation (\ref{p.3}) the space stops to be Euclidean and becomes
a Riemannian space $R_n=\left\{ \mathbf{\tilde g},\tilde K,\Bbb{R}^n\right\}
$. Constraint (\ref{p.3}) is a condition of the Euclideness of the space.

Eliminating (\ref{p.3}) one obtains a Riemannian space $R_n=\left\{ \mathbf{g%
},K,\Bbb{R}^n\right\} $, which is determined by the form of the metric
tensor $g_{ik}(x)$. The world function is determined by the relation (\ref
{a3.15}), where $\mathcal{L}_{[xx^{\prime }]}\subset \Bbb{R}^n$ is the
geodesic segment of the geodesic $\mathcal{L}_{xx^{\prime }}\subset \Bbb{R}%
^n $. This geodesic is an extremal of (\ref{a3.15}), considered as a
functional of the curve $\mathcal{L}:\;x=x(\tau )$, written in the form
\begin{equation}
\sigma [x(\tau )]=\frac 12\left( \int\limits_{_{\mathcal{L}}}\sqrt{g_{ik}(x)%
\frac{dx^i}{d\tau }\frac{dx^k}{d\tau }}d\tau \right) ^2  \label{a3.15a}
\end{equation}
The geodesic $\mathcal{L}_{xx^{\prime }}:\;x=x(\tau )$ is described by the
equations
\begin{equation}
\mathcal{L}_{xx^{\prime }}:\qquad \frac{d^2x^i}{d\tau ^2}+\gamma _{kl}^i%
\frac{dx^k}{d\tau }\frac{dx^l}{d\tau }=0,\qquad i=1,2,\ldots n  \label{a.2.2}
\end{equation}
where
\begin{equation}
\gamma _{kl}^i=\gamma _{kl}^i(x)=\frac 12g^{ij}\left(
g_{kj,l}+g_{lj,k}-g_{kl,j}\right)  \label{a.2.3}
\end{equation}
is the Christoffel symbol, and comma before index $l$ denotes
differentiation with respect to $x^l$.

In particular, if $g_{ik}=$const, $i,k=1,2,\ldots n$, $\;g=\det
||g_{ik}||\neq 0$, the world function is described by the relation (\ref
{a3.14}), and the Riemannian space $R_n=\left\{ \mathbf{g},K,\Bbb{R}%
^n\right\} $ is the Euclidean space. Let us consider now the Riemannian
space $R_n=\left\{ \mathbf{g}_E,K,D\right\} $, where $D\subset \Bbb{R}^n$ is
some region of the Euclidean space $E_n=\left\{ \mathbf{g}_E,K,\Bbb{R}%
^n\right\} .$ If the region $D$ is convex, i.e. any segment $\mathcal{L}%
_{[xx^{\prime }]}$ of the straight $\mathcal{L}_{xx^{\prime }}$, connecting
the points $x,x^{\prime }\in D$ belongs to $D$ ($\mathcal{L}_{[xx^{\prime
}]}\subset D)$, then the world function of the Riemannian space $R_n=\left\{
\mathbf{g}_E,K,D\right\} $ has the form (\ref{a3.14}) and the Riemannian
space $R_n=\left\{ \mathbf{g}_E,K,D\right\} $ can be embedded isometrically
into the Euclidean space $E_n=\left\{ \mathbf{g}_E,K,\Bbb{R}^n\right\} $.

If the region $D$ is noncovex, then the system of geodesics of $R_n=\left\{
\mathbf{g}_E,K,D\right\} $ is not a system of straight lines, and the world
function (\ref{a3.15}) is not described by the relation (\ref{a3.14}).

\textit{Example}. Let us consider two-dimensional proper Euclidean space,
and rectilinear orthogonal coordinates on it. Let us consider the region $%
D:\;$ $\left( x^1\right) ^2+\left( x^2\right) ^2\geq 1.$ Geodesics of the
Riemannian space $R_2^{\prime }=\left\{ \mathbf{g}_E,K,D\right\} $ looks as
it is shown in Figure 1. After cutting a hole in the Euclidean plane the
shape and length of geodesic segment between the points $P$ and $P^{\prime }$
changes. World function $\sigma (P,P^{\prime })$ between the points $P$ and $%
P^{\prime }$ changes, and the part $R_2^{\prime }=\left\{ \mathbf{g}%
_E,K,D\right\} $ of the Euclidean plane $R_2=\left\{ \mathbf{g}_E,K,\Bbb{R}%
^2\right\} $ stops to be embeddable isometrically in $R_2=\left\{ \mathbf{g}%
_E,K,\Bbb{R}^2\right\} $. It seems to be rather strange, when part of the
Euclidean plane cannot be embedded isometrically in the plane.

The problem of convexity is rather strong, and most of geometricians prefer
to get around this problem, considering convex regions \cite{A48}. In the
T-geometry the convexity problem is absent. Indeed, according to definition
\ref{d3.1.1ae} any subset of a $\sigma $-space is always embeddable
isometrically into the $\sigma $-space. From viewpoint of T-geometry,
cutting a hole in the Euclidean plane $R_2=\left\{ \mathbf{g}_E,K,\Bbb{R}%
^2\right\} $, one does not change the system of geodesics (the first order
NGOs), one cuts only holes in geodesics, making them discontinuous.
Continuity is a property of coordinate systems, used in Riemannian geometry
as the main tool of description. From viewpoint of T-geometry the convexity
problem is a problem made artificially. Insisting on continuity of
geodesics, one overestimates importance of the continuity for geometry and
attributes the continuity of geodesics (the first order NGOs) to any
Riemannian geometry, whereas the continuity of geodesics is a special
property of the proper Euclidean geometry.

\section{Riemannian geometry and one-dimensionality of the first order NGOs}

Let us consider the $n$-dimensional pseudo-Euclidean space $E_n=\left\{
\mathbf{g}_1,K,\Bbb{R}^n\right\} $ of the index $1$, $\mathbf{g}_1=$diag$%
\left\{ 1,-1,-1\ldots -1\right\} $ to be a kind of $n$-dimensional
Riemannian space\footnote{%
The term ''Riemannian space'' is considered to be a collective term with
respect to concepts ''proper Riemannian'' and ''pseudo-Riemannian''. Matrix $%
\mathbf{g}$ of the metric tensor has eigenvalues of the same sign in the
case of proper Riemannian space and of different signs in the case of
pseudo-Riemannian one.}. The world function is defined by the relation (\ref
{a3.14})
\begin{equation}
\sigma _1(x,x^{\prime })={\frac 12}\sum_{i,k=1}^ng^{ik}\left(
x_i-x_i^{\prime }\right) \left( x_k-x_k^{\prime }\right) ,\qquad g^{ik}=%
\mathrm{diag}\left\{ 1,-1,-1\ldots -1\right\}  \label{c3.0}
\end{equation}
Geodesic $\mathcal{L}_{yy^{\prime }}$ is a straight line, and it is
considered in pseudo-Euclidean geometry to be the first order NGOs,
determined by two points $y$ and $y^{\prime }$

\begin{equation}
\mathcal{L}_{yy^{\prime }}:\quad x^i=\left( y^i-y^{\prime i}\right) \tau
,\qquad i=1,2,\ldots n,\qquad \tau \in \Bbb{R}  \label{c3.1}
\end{equation}
The geodesic $\mathcal{L}_{yy^{\prime }}$ is called timelike, if $\sigma
_1(y,y^{\prime })>0$, and it is called spacelike if $\sigma _1(y,y^{\prime
})<0$. The geodesic $\mathcal{L}_{yy^{\prime }}$ is called null, if $\sigma
_1(y,y^{\prime })=0$.

The pseudo-Euclidean space $E_n=\left\{ \mathbf{g}_1,K,\Bbb{R}^n\right\} $
generates the $\sigma $-space $V=\left\{ \sigma _1,\Bbb{R}^n\right\} $,
where the world function $\sigma _1$ is defined by the relation (\ref{c3.0}%
). The first order tube (NGO) $\mathcal{T}\left( x,x^{\prime }\right) $ in
the $\sigma $-Riemannian space $V=\left\{ \sigma _1,\Bbb{R}^n\right\} $ is
defined by the relation (\ref{b1.3})
\begin{equation}
\mathcal{T}\left( x,x^{\prime }\right) \equiv \mathcal{T}_{xx^{\prime
}}=\left\{ r|F_2\left( x,x^{\prime },r\right) =0\right\} ,\qquad \sigma
_1(x,x^{\prime })\neq 0,\qquad x,x^{\prime },r\in \Bbb{R}^n,  \label{c3.2}
\end{equation}
\begin{equation}
F_2\left( x,x^{\prime },r\right) =\left|
\begin{array}{cc}
(x_i^{\prime }-x_i)(x^{\prime i}-x^i) & (x_i^{\prime }-x_i)(r^i-x^i) \\
(r_i-x_i)(x^{\prime i}-x^i) & (r_i-x_i)(r^i-x^i)
\end{array}
\right|  \label{c3.3}
\end{equation}
Solution of equations (\ref{c3.2}), (\ref{c3.3}) gives the following result
\begin{equation}
\mathcal{T}_{xx^{\prime }}=\left\{ r\left| \bigcup\limits_{y\in \Bbb{R}%
^n}\bigcup\limits_{\tau \in \Bbb{R}}r=\left( x^{\prime }-x\right) \tau
+y-x\wedge \Gamma (x,x^{\prime },y)=0\right. \wedge \Gamma (x,y,y)=0\right\}
,  \label{c3.4}
\end{equation}
\[
x,x^{\prime },y,r\in \Bbb{R}^n
\]
where $\Gamma (x,x^{\prime },y)=(x_i^{\prime }-x_i)(y^i-x^i)$ is the scalar
product of vectors $\overrightarrow{xy}$ and $\overrightarrow{xx^{\prime }}$
defined by the relation (\ref{a1.7}). In the case of timelike vector $%
\overrightarrow{xx^{\prime }}$, when $\sigma _1(x,x^{\prime })>0$, there is
a unique null vector $\overrightarrow{xy}=\overrightarrow{xx}=%
\overrightarrow{0}$ which is orthogonal to the vector $\overrightarrow{%
xx^{\prime }}$. In this case the ($n-1)$-dimensional surface $\mathcal{T}%
_{xx^{\prime }}$ degenerates into the one-dimensional straight
\begin{equation}
\mathcal{T}_{xx^{\prime }}=\left\{ r\left| \bigcup\limits_{\tau \in \Bbb{R}%
}r=\left( x^{\prime }-x\right) \tau \right. \right\} ,\qquad \sigma
_1(x,x^{\prime })>0,\qquad x,x^{\prime },r\in \Bbb{R}^n,  \label{c3.5}
\end{equation}
Thus, for timelike vector $\overrightarrow{xx^{\prime }}$ the first order
tube $\mathcal{T}_{xx^{\prime }}$ coincides with the geodesic $\mathcal{L}%
_{xx^{\prime }}$. In the case of spacelike vector $\overrightarrow{%
xx^{\prime }}$ the $(n-1)$-dimensional tube $\mathcal{T}_{xx^{\prime }}$
contains the one-dimensional geodesic $\mathcal{L}_{xx^{\prime }}$ of the
pseudo-Euclidean space $E_n=\left\{ \mathbf{g}_1,K,\Bbb{R}^n\right\} $.

This difference poses the question what is the reason of this difference and
what of the two generalization of the proper Euclidean geometry is more
reasonable. Note that four-dimensional pseudo-Euclidean geometry is used for
description of the real space-time. One can try to resolve this problem from
experimental viewpoint. Free classical particles are described by means of
timelike straight lines. At this point the pseudo-Euclidean geometry and the
$\sigma $-pseudo-Euclidean geometry (T-geometry) lead to the same result.
The spacelike straights are believed to describe the particles moving with
superlight speed (so-called taxyons). Experimental attempts of taxyons
discovery were failed. Of course, trying to discover taxyons, one considered
them to be described by spacelike straights. On the other hand, the
physicists believe that all what can exist does exist and may be discovered.
From this viewpoint the failure of discovery of taxyons in the form of
spacelike line justifies in favour of taxyons in the form of
three-dimensional surfaces.

To interpret the structure of the set (\ref{c3.4}), describing the first
order tube, let us take into account the zeroth order tube $\mathcal{T}_x$,
determined by the point $x$ in the $\sigma $-pseudo-Euclidean space is the
light cone with the vertex at the point $x$ (not the point $x$). Practically
the first order tube consists of such sections of the light cones with their
vertex $y\in \mathcal{L}_{xx^{\prime }}$ that all vectors $\overrightarrow{yr%
}$ of these sections are orthogonal to the vector $\overrightarrow{%
xx^{\prime }} $. In other words, the first order tube $\mathcal{T}%
_{xx^{\prime }}$ consists of the zeroth order tubes $\mathcal{T}_y$ sections
at $y$, orthogonal to $\overrightarrow{xx^{\prime }}$, with $y\in \mathcal{L}%
_{xx^{\prime }}$. For timelike $\overrightarrow{xx^{\prime }}$ this section
consists of one point, but for the spacelike $\overrightarrow{xx^{\prime }}$
it is two-dimensional section of the light cone.

\section{Collinearity in Riemannian and $\sigma $-Riemannian geometry}

Let us return to the Riemannian space $R_n=\left\{ \mathbf{g},K,D\right\}
,\,\quad D\subset \Bbb{R}^n$, which generates the world function $\sigma
(x,x^{\prime })$ defined by the relation (\ref{a3.15}). Then the $\sigma $%
-space $V=\left\{ \sigma ,D\right\} $ appears. it will be referred to as $%
\sigma $-Riemannian space. We are going to compare concept of collinearity
(parallelism) of two vectors in the two spaces.

The world function $\sigma =\sigma (x,x^{\prime })$ of both $\sigma $%
-Riemannian and Riemannian spaces satisfies the system of equations \cite
{R62}\footnote{%
The paper \cite{R62} is hardly available for English speaking reader. Survey
of main results of \cite{R62} in English may be found in \cite{R64}. See
also \cite{R92}}
\begin{equation}
\begin{array}{cc}
(1)\quad \sigma _l\sigma ^{{lj^{\prime }}}\sigma _{j^{\prime }}=2\sigma
\quad \quad \quad & (4)\quad \det \parallel \sigma _{i||k}\parallel \neq 0
\\
(2)\quad \sigma (x,x^{\prime })=\sigma (x^{\prime },x)\quad \quad & (5)\quad
\det \parallel \sigma _{{ik^{\prime }}}\parallel \neq 0 \\
(3)\quad \sigma (x,x)=0\quad \quad \quad \quad \quad & (6)\quad \sigma
_{i||k||l}=0\quad \quad
\end{array}
\label{b2.39}
\end{equation}
where the following designations are used
\[
\sigma _i\equiv \frac{\partial \sigma }{\partial x^i},\qquad \sigma
_{i^{\prime }}\equiv \frac{\partial \sigma }{\partial x^{\prime i}},\qquad
\sigma _{ik^{\prime }}\equiv \frac{\partial ^2\sigma }{\partial x^i\partial
x^{\prime k}},\qquad \sigma ^{ik^{\prime }}\sigma _{lk^{\prime }}=\delta
_l^i
\]
Here the primed index corresponds to the point $x,$ and unprimed index
corresponds to the point $x$. Two parallel vertical strokes mean covariant
derivative $\tilde \nabla _i^{x^{\prime }}$ with respect to $x^i$ with the
Christoffel symbol
\[
\Gamma _{kl}^i\equiv \Gamma _{kl}^i\left( x,x^{\prime }\right) \equiv \sigma
^{is^{\prime }}\sigma _{kls^{\prime }},\qquad \sigma _{kls^{\prime }}\equiv
\frac{\partial ^3\sigma }{\partial x^k\partial x^l\partial x^{\prime s}}
\]
For instance,
\begin{equation}
G_{ik}\equiv G_{ik}(x,x^{\prime })\equiv \sigma _{i||k}\equiv \frac{\partial
\sigma _i}{\partial x^k}-\Gamma _{ik}^l\left( x,x^{\prime }\right) \sigma
_l\equiv \frac{\partial \sigma _i}{\partial x^k}-\sigma _{iks^{\prime
}}\sigma ^{ls^{\prime }}\sigma _l  \label{a2.0}
\end{equation}
\[
G_{ik||l}\equiv \frac{\partial G_{ik}}{\partial x^l}-\sigma _{ils^{\prime
}}\sigma ^{js^{\prime }}G_{jk}-\sigma _{kls^{\prime }}\sigma ^{js^{\prime
}}G_{ij}
\]
Summation from $1$ to $n$ is produced over repeated indices. The covariant
derivative $\tilde \nabla _i^{x^{\prime }}$ with respect to $x^i$ with the
Christoffel symbol $\Gamma _{kl}^i\left( x,x^{\prime }\right) $ acts only on
the point $x$ and on unprimed indices. It is called the tangent derivative,
because it is a covariant derivative in the Euclidean space $E_{x^{\prime }}$
which is tangent to the Riemannian space $R_n$ at the point $x^{\prime }$.
The covariant derivative $\tilde \nabla _{i^{\prime }}^x$ with respect to $%
x^{\prime i}$ with the Christoffel symbol $\Gamma _{k^{\prime }l^{\prime
}}^{i^{\prime }}\left( x,x^{\prime }\right) $ acts only on the point $%
x^{\prime }$ and on primed indices. It is a covariant derivative in the
Euclidean space $E_x$ which is tangent to the $\sigma $-Riemannian space $V$
at the point $x$ \cite{R62}.

In general, the world function $\sigma $ carries out the geodesic mapping $%
G_{x^{\prime }}:$ $R_n\rightarrow E_{x^{\prime }}$ of the Riemannian space $%
R_n=\left\{ \mathbf{g},K,D\right\} $ on the Euclidean space $E_{x^{\prime
}}=\left\{ \mathbf{g},K_{x^{\prime }},D\right\} $, tangent to $R_n=\left\{
\mathbf{g},K,D\right\} $ at the point $x^{\prime }$ \cite{R62}. This mapping
transforms the coordinate system $K$ in $R_n$ into the coordinate system $%
K_{x^{\prime }}$ in $E_{x^{\prime }}$. The mapping is geodesic in the sense
that it conserves the lengths of segments of all geodesics, passing through
the tangent point $x^{\prime }$ and angles between them at this point.

The tensor $G_{ik}$, defined by (\ref{a2.0}) is the metric tensor at the
point $x$ in the tangent Euclidean space $E_{x^{\prime }}$. The covariant
derivatives $\tilde \nabla _i^{x^{\prime }}$ and $\tilde \nabla
_k^{x^{\prime }}$ commute identically, i.e. $(\tilde \nabla _i^{x^{\prime }}
\tilde \nabla _k^{x^{\prime }}-\tilde \nabla _k^{x^{\prime }}\tilde \nabla
_i^{x^{\prime }})A_{ls}\equiv 0,$ for any tensor $A_{ls}$ \cite{R62}. This
shows that they are covariant derivatives in the flat space $E_{x^{\prime }}$%
.

The system of equations (\ref{b2.39}) contains only world function $\sigma $
and its derivatives, nevertheless the system of equations (\ref{b2.39}) is
not $\sigma $-immanent, because it contains a reference to a coordinate
system. It does not contain the metric tensor explicitly. Hence, it is valid
for any Riemannian space $R_n=\left\{ \mathbf{g},K,D\right\} $. All
relations written above are valid also for the $\sigma $-space $V=\left\{
\sigma ,D\right\} $, provided the world function $\sigma $ is coupled with
the metric tensor by relation (\ref{a3.15}).

$\sigma $-immanent expression for scalar product $\left( \mathbf{P}_0\mathbf{%
P}_1.\mathbf{Q}_0\mathbf{Q}_1\right) $ of two vectors $\mathbf{P}_0\mathbf{P}%
_1$ and $\mathbf{Q}_0\mathbf{Q}_1$ in the proper Euclidean space has the
form
\begin{equation}
\left( \mathbf{P}_0\mathbf{P}_1.\mathbf{Q}_0\mathbf{Q}_1\right) \equiv
\sigma \left( P_0,Q_1\right) +\sigma \left( Q_0,P_1\right) -\sigma \left(
P_0,Q_0\right) -\sigma \left( P_1,Q_1\right)  \label{a2.1a}
\end{equation}
This relation can be easily proved as follows.

In the proper Euclidean space three vectors $\mathbf{P}_0\mathbf{P}_1$, $%
\mathbf{P}_0\mathbf{Q}_1$, and $\mathbf{P}_1\mathbf{Q}_1$ are coupled by the
relation
\begin{equation}
\mid \mathbf{P}_1\mathbf{Q}_1\mid ^2=\mid \mathbf{P}_0\mathbf{Q}_1-\mathbf{P}%
_0\mathbf{P}_1\mid ^2=\mid \mathbf{P}_0\mathbf{P}_1\mid ^2+\mid \mathbf{P}_0%
\mathbf{Q}_1\mid ^2-2(\mathbf{P}_0\mathbf{P}_1.\mathbf{P}_0\mathbf{Q}_1)
\label{a2.15}
\end{equation}
where $(\mathbf{P}_0\mathbf{P}_1.\mathbf{P}_0\mathbf{Q}_1)$ denotes the
scalar product of two vectors $\mathbf{P}_0\mathbf{P}_1$ and $\mathbf{P}_0%
\mathbf{Q}_1$ in the proper Euclidean space. It follows from (\ref{a2.15})
\begin{equation}
(\mathbf{P}_0\mathbf{P}_1.\mathbf{P}_0\mathbf{Q}_1)={\frac 12}\{\mid \mathbf{%
P}_0\mathbf{Q}_1\mid ^2+\mid \mathbf{P}_0\mathbf{P}_1\mid ^2-\mid \mathbf{P}%
_1\mathbf{Q}_1\mid ^2\}  \label{a2.16}
\end{equation}
Substituting the point $Q_1$ by $Q_0$ in (\ref{a2.16}), one obtains
\begin{equation}
(\mathbf{P}_0\mathbf{P}_1.\mathbf{P}_0\mathbf{Q}_0)={\frac 12}\{\mid \mathbf{%
P}_0\mathbf{Q}_0\mid ^2+\mid \mathbf{P}_0\mathbf{P}_1\mid ^2-\mid \mathbf{P}%
_1\mathbf{Q}_0\mid ^2\}  \label{a2.16a}
\end{equation}
Subtracting (\ref{a2.16a}) from (\ref{a2.16}) and using the properties of
the scalar product in the Euclidean space, one obtains
\begin{equation}
(\mathbf{P}_0\mathbf{P}_1.\mathbf{Q}_0\mathbf{Q}_1)={\frac 12}\{\mid \mathbf{%
P}_0\mathbf{Q}_1\mid ^2+\mid \mathbf{Q}_0\mathbf{P}_1\mid ^2-\mid \mathbf{P}%
_0\mathbf{Q}_0\mid ^2-\mid \mathbf{P}_1\mathbf{Q}_1\mid ^2\}  \label{a2.16b}
\end{equation}
Taking into account that $\mid \mathbf{P}_0\mathbf{Q}_1\mid ^2=2\sigma
\left( P_0,Q_1\right) $, one obtains the relation\ (\ref{a2.1a}) from the
relation (\ref{a2.16b}).

Two vectors $\mathbf{P}_0\mathbf{P}_1$ and $\mathbf{Q}_0\mathbf{Q}_1$ are
collinear $\mathbf{P}_0\mathbf{P}_1||\mathbf{Q}_0\mathbf{Q}_1$ (parallel or
antiparallel), provided $\cos ^2\theta =1,$ where $\theta $ is the angle
between the vectors $\mathbf{P}_0\mathbf{P}_1$ and $\mathbf{Q}_0\mathbf{Q}_1$%
. Taking into account that

\begin{equation}
\cos ^2\theta =\frac{\left( \mathbf{P}_0\mathbf{P}_1.\mathbf{Q}_0\mathbf{Q}%
_1\right) ^2}{\left( \mathbf{P}_0\mathbf{P}_1.\mathbf{P}_0\mathbf{P}%
_1\right) (\mathbf{Q}_0\mathbf{Q}_1.\mathbf{Q}_0\mathbf{Q}_1)}=\frac{\left(
\mathbf{P}_0\mathbf{P}_1.\mathbf{Q}_0\mathbf{Q}_1\right) ^2}{|\mathbf{P}_0%
\mathbf{P}_1|^2\cdot |\mathbf{Q}_0\mathbf{Q}_1|^2}  \label{a2.2}
\end{equation}
one obtains the following $\sigma $-immanent condition of the two vectors
collinearity
\begin{equation}
\mathbf{P}_0\mathbf{P}_1||\mathbf{Q}_0\mathbf{Q}_1:\qquad \left( \mathbf{P}_0%
\mathbf{P}_1.\mathbf{Q}_0\mathbf{Q}_1\right) ^2=|\mathbf{P}_0\mathbf{P}%
_1|^2\cdot |\mathbf{Q}_0\mathbf{Q}_1|^2  \label{a2.3}
\end{equation}
The collinearity condition (\ref{a2.3}) is $\sigma $-immanent, because by
means of (\ref{a2.1a}) it can be written in terms of the $\sigma $-function
only. Thus, this relation describes the vectors collinearity in the case of
arbitrary $\sigma $-space.

Let us describe this relation for the case of $\sigma $-Riemannian geometry.
Let coordinates of the points $P_0,P_1,Q_0,Q_1$ be respectively $x,$ $x+dx,$
$x^{\prime }$ and $x^{\prime }+dx^{\prime }$. Then writing (\ref{a2.1a}) and
expanding it over $dx$ and $dx^{\prime }$, one obtains
\begin{eqnarray*}
\left( \mathbf{P}_0\mathbf{P}_1.\mathbf{Q}_0\mathbf{Q}_1\right) &\equiv
&\sigma \left( x,x^{\prime }+dx^{\prime }\right) +\sigma \left( x^{\prime
},x+dx\right) -\sigma \left( x,x^{\prime }\right) -\sigma \left(
x+dx,x^{\prime }+dx^{\prime }\right) =  \label{a2.1b} \\
&&\ \sigma +\sigma _{l^{\prime }}dx^{\prime l^{\prime }}+\frac 12\sigma
_{l^{\prime },s^{\prime }}dx^{\prime l^{\prime }}dx^{\prime s^{\prime
}}+\sigma +\sigma _idx^i+\frac 12\sigma _{i,k}dx^idx^k-\sigma \\
&&\ -\sigma -\sigma _idx^i-\sigma _{l^{\prime }}dx^{\prime l^{\prime
}}-\frac 12\sigma _{i,k}dx^idx^k-\sigma _{i,l^{\prime }}dx^idx^{\prime
l^{\prime }}-\frac 12\sigma _{l^{\prime },s^{\prime }}dx^{\prime l^{\prime
}}dx^{\prime s^{\prime }}
\end{eqnarray*}
\begin{equation}
\left( \mathbf{P}_0\mathbf{P}_1.\mathbf{Q}_0\mathbf{Q}_1\right) =-\sigma
_{i,l^{\prime }}dx^idx^{\prime l^{\prime }}=-\sigma _{il^{\prime
}}dx^idx^{\prime l^{\prime }}  \label{a2.4}
\end{equation}
Here comma means differentiation. For instance, $\sigma _{i,k}\equiv
\partial \sigma _i/\partial x^k$. One obtains for $|\mathbf{P}_0\mathbf{P}%
_1|^2$ and $|\mathbf{Q}_0\mathbf{Q}_1|^2$%
\begin{equation}
|\mathbf{P}_0\mathbf{P}_1|^2=g_{ik}dx^idx^k,\qquad |\mathbf{Q}_0\mathbf{Q}%
_1|^2=g_{l^{\prime }s^{\prime }}dx^{\prime l^{\prime }}dx^{\prime s^{\prime
}}  \label{a2.5}
\end{equation}
where $g_{ik}=g_{ik}(x)$ and $g_{l^{\prime }s^{\prime }}=g_{l^{\prime
}s^{\prime }}(x^{\prime })$. Then the collinearity condition (\ref{a2.3}) is
written in the form
\begin{equation}
\left( \sigma _{il^{\prime }}\sigma _{ks^{\prime }}-g_{ik}g_{l^{\prime
}s^{\prime }}\right) dx^idx^kdx^{\prime l^{\prime }}dx^{\prime s^{\prime }}=0
\label{a4.3}
\end{equation}
Let us take into account that in the Riemannian space the metric tensor $%
g_{l^{\prime }s^{\prime }}$ at the point $x^{\prime }$ can be expressed via
the world function $\sigma $ of points $x,x^{\prime }$ by means of the
relation \cite{R62}
\begin{equation}
g_{l^{\prime }s^{\prime }}=\sigma _{il^{\prime }}G^{ik}\sigma _{ks^{\prime
}},\qquad g^{l^{\prime }s^{\prime }}=\sigma ^{il^{\prime }}G_{ik}\sigma
^{ks^{\prime }}  \label{a4.4}
\end{equation}
where the tensor $G_{ik}$ is defined by the relation (\ref{a2.0}), and $%
G^{ik}$ is defined by the relation
\begin{equation}
G^{il}G_{lk}=\delta _k^i  \label{a4.5}
\end{equation}

Substituting the first relation (\ref{a2.0}) in (\ref{a4.3}) and using
designation
\begin{equation}
u_i=-\sigma _{il^{\prime }}dx^{\prime l^{\prime }},\qquad
u^i=G^{ik}u_k=-\sigma ^{il^{\prime }}g_{l^{\prime }s^{\prime }}dx^{l^{\prime
}s^{\prime }}  \label{a4.6}
\end{equation}
one obtains
\begin{equation}
\left( \delta _i^l\delta _k^s-g_{ik}G^{ls}\right) u_lu_sdx^idx^k=0
\label{a2.6}
\end{equation}
The vector $u_i$ is the vector $dx_{i^{\prime }}^{\prime }=g_{i^{\prime
}k^{\prime }}dx^{\prime k^{\prime }}$ transported parallelly from the point $%
x^{\prime }$ to the point $x$ in the Euclidean space $E_{x^{\prime }}$
tangent to the Riemannian space $R_n$. Indeed,
\begin{equation}
u_i=-\sigma _{il^{\prime }}g^{l^{\prime }s^{\prime }}dx_s^{\prime },\qquad
\tilde \nabla _k^{x^{\prime }}\left(-\sigma _{il^{\prime }}g^{l^{\prime
}s^{\prime }}\right) \equiv 0,\qquad i,k=1,2,\ldots n  \label{a2.7}
\end{equation}
and tensor $-\sigma _{il^{\prime }}g^{l^{\prime }s^{\prime }}$ is the
operator of the parallel transport in $E_{x^{\prime }}$, because
\[
\left[-\sigma _{il^{\prime }}g^{l^{\prime }s^{\prime
}}\right]_{x=x^{\prime}}= \delta ^{s^{\prime}}_{i^{\prime}}
\]
and the tangent derivative of this operator is equal to zero identically.
For the same reason, i.e. because of
\[
\left[\sigma ^{il^{\prime }}g_{l^{\prime }s^{\prime }}\sigma ^{ks^{\prime
}}\right]_{x=x^{\prime}}=g^{i^{\prime}k^{\prime}}, \qquad \tilde \nabla
_s^{x^{\prime }}(\sigma ^{il^{\prime }}g_{l^{\prime }s^{\prime }}\sigma
^{ks^{\prime }})\equiv 0
\]
$G^{ik}=\sigma ^{il^{\prime }}g_{l^{\prime }s^{\prime }}\sigma ^{ks^{\prime
}}$ is the contravariant metric tensor in $E_{x^{\prime }}$, at the point $x
$.

The relation (\ref{a2.6}) contains vectors at the point $x$ only . At fixed $%
u_i=-\sigma _{il^{\prime }}dx^{\prime l^{\prime }}$ it describes a
collinearity cone, i.e. a cone of infinitesimal vectors $dx^i$ at the point $%
x $ parallel to the vector $dx^{\prime i^{\prime }}$ at the point $x^{\prime
}$. Under some condition the collinearity cone can degenerates into a line.
In this case there is only one direction, parallel to the fixed vector $u^i$%
. Let us investigate, when this situation takes place.

At the point $x$ two metric tensors $g_{ik}$ and $G_{ik}$ are connected by
the relation \cite{R62}
\begin{equation}
G_{ik}(x,x^{\prime })=g_{ik}(x)+\int\limits_x^{x^{\prime }}F_{{ikj^{\prime
\prime }}s^{{\prime \prime }}}(x,x^{\prime \prime })\sigma ^{j^{{\prime
\prime }}}(x,x^{\prime \prime })d{x^{\prime \prime }}^{s^{\prime \prime }},
\label{b3.23}
\end{equation}
where according to \cite{R62}
\begin{equation}
\sigma ^{i^{\prime }}=\sigma ^{li^{\prime }}\sigma _l=G^{l^{\prime
}i^{\prime }}\sigma _{l^{\prime }}=g^{l^{\prime }i^{\prime }}\sigma
_{l^{\prime }}  \label{a4.6a}
\end{equation}
Integration does not depend on the path, because it is produced in the
Euclidean space $E_{x^{\prime }}$. The two-point tensor $F_{{ilk^{\prime
}j^{\prime }}}=F_{{ilk^{\prime }j^{\prime }}}(x,x^{\prime })$ is the
two-point curvature tensor, defined by the relation
\begin{equation}
F_{{ilk^{\prime }j^{\prime }}}=\sigma _{{ilj^{\prime }}\parallel k^{\prime
}}=\sigma _{{ilj^{\prime }},k^{\prime }}-\sigma _{{sj^{\prime }k^{\prime }}%
}\sigma ^{{sm^{\prime }}}\sigma _{{ilm^{\prime }}}=\sigma _{i\mid
l||k^{\prime }||j^{\prime }}  \label{b3.5}
\end{equation}
where one vertical stroke denotes usual covariant derivative and two
vertical strokes denote tangent derivative. The two-point curvature tensor $%
F_{{ilk^{\prime }j^{\prime }}}$ has the following symmetry properties
\begin{equation}
F_{{ilk^{\prime }j^{\prime }}}=F_{{lik^{\prime }j^{\prime }}}=F_{{%
ilj^{\prime }k^{\prime }}},\qquad F_{{ilk^{\prime }j^{\prime }}}(x,x^{\prime
})=F_{{k^{\prime }j^{\prime }il}}(x^{\prime },x)  \label{b3.10}
\end{equation}
It is connected with the one-point Riemann-Ghristoffel curvature tensor $%
r_{iljk}$ by means of relations
\begin{equation}
r_{{iljk}}=\left[ F_{{ikj^{\prime }l^{\prime }}}-F_{{ijk^{\prime }l^{\prime }%
}}\right] _{x^{\prime }=x}=f_{{ikjl}}-f_{{ijkl}},\qquad f_{{iklj}}=\left[ F_{%
{ikj^{\prime }l^{\prime }}}\right] _{x^{\prime }=x}  \label{b3.14}
\end{equation}

In the Euclidean space the two-point curvature tensor $F_{{ilk^{\prime
}j^{\prime }}}$ vanishes as well as the Riemann-Ghristoffel curvature tensor
$r_{iljk}$.

Let us introduce designation
\begin{equation}
\Delta _{ik}=\Delta _{ik}(x,x^{\prime })=\int\limits_x^{x^{\prime }}F_{{%
ikj^{\prime \prime }}s^{{\prime \prime }}}(x,x^{\prime \prime })\sigma ^{j^{{%
\prime \prime }}}(x,x^{\prime \prime })d{x^{\prime \prime }}^{s^{\prime
\prime }}  \label{a4.8}
\end{equation}
and choose the geodesic $\mathcal{L}_{xx^{\prime }}$ as the path of
integration. It is described by the relation
\begin{equation}
\sigma _i(x,x^{\prime \prime })=\tau \sigma _i(x,x^{\prime })  \label{a4.11}
\end{equation}
which determines $x^{\prime \prime }$ as a function of parameter $\tau $.
Differentiating with respect to $\tau $, one obtains
\begin{equation}
\sigma _{ik^{\prime \prime }}(x,x^{\prime \prime })dx^{\prime \prime
k^{\prime \prime }}=\sigma _i(x,x^{\prime })d\tau  \label{a4.10}
\end{equation}
Resolving equations (\ref{a4.10}) with respect to $dx^{\prime \prime }$ and
substituting in (\ref{a4.8}), one obtains
\begin{equation}
\Delta _{ik}(x,x^{\prime })=\sigma _l(x,x^{\prime })\sigma _p(x,x^{\prime
})\int\limits_0^1F_{{ikj^{\prime \prime }}s^{{\prime \prime }}}(x,x^{\prime
\prime })\sigma ^{lj^{{\prime \prime }}}(x,x^{\prime \prime })\sigma ^{ps^{{%
\prime \prime }}}(x,x^{\prime \prime })\tau d{\tau }  \label{a4.12}
\end{equation}
where $x^{\prime \prime }$ is determined from (\ref{a4.11}) as a function of
$\tau $. Let us set
\begin{equation}
F_{{ik}}^{\;..lp}(x,x^{\prime })=F_{{ikj^{\prime }}s^{{\prime }%
}}(x,x^{\prime })\sigma ^{lj^{{\prime }}}(x,x^{\prime })\sigma ^{ps^{{\prime
}}}(x,x^{\prime })  \label{a4.13}
\end{equation}
then
\begin{equation}
G_{ik}(x,x^{\prime })=g_{ik}(x)+\Delta _{ik}(x,x^{\prime })  \label{a4.14}
\end{equation}
\begin{equation}
\Delta _{ik}(x,x^{\prime })=\sigma _l(x,x^{\prime })\sigma _p(x,x^{\prime
})\int\limits_0^1F_{{ik}}^{\;..lp}(x,x^{\prime \prime })\tau d{\tau }
\label{a4.15}
\end{equation}
Substituting $g_{ik}$ from (\ref{a4.14}) in (\ref{a2.6}), one obtains
\begin{equation}
\left( \delta _i^l\delta _k^s-G^{ls}\left( G_{ik}-\Delta _{ik}\right)
\right) u_lu_sdx^idx^k=0  \label{a4.7}
\end{equation}

Let us look for solutions of equation in the form of expansion
\begin{equation}
dx^i=\alpha u^i+v^i,\qquad G_{ik}u^iv^k=0  \label{a4.16}
\end{equation}
Substituting (\ref{a4.16}) in (\ref{a4.7}), one obtains equation for $v^i$%
\begin{equation}
G_{ls}u^lu^s\left[ G_{ik}v^iv^k-\Delta _{ik}\left( \alpha u^i+v^i\right)
\left( \alpha u^k+v^k\right) \right] =0  \label{a4.18}
\end{equation}
If the $\sigma $-Riemannian space $V=\left\{ \sigma ,D\right\} $ is $\sigma $%
-Euclidean, then as it follows from (\ref{a4.15}) $\Delta _{ik}=0$. If $%
V=\left\{ \sigma ,D\right\} $ is the proper $\sigma $-Euclidean space, $%
G_{ls}u^lu^s\neq 0$, and one obtains two equations for determination of $v^i$
\begin{equation}
G_{ik}v^iv^k=0,\qquad G_{ik}u^iv^k=0  \label{a4.17}
\end{equation}
The only solution
\begin{equation}
v^i=0,\qquad dx^i=\alpha u^i,\qquad i=1,2,\ldots n  \label{a2.6b}
\end{equation}
of (\ref{a4.18}) is a solution of the equation (\ref{a4.7}), where $\alpha $
is an arbitrary constant. In the proper Euclidean geometry the collinearity
cone always degenerates into a line.

Let now the space $V=\left\{ \sigma ,D\right\} $ be the $\sigma $%
-pseudo-Euclidean space of index $1$, and the vector $u^i$ be timelike, i.e.
$G_{ik}u^iu^k>0$. Then equations (\ref{a4.17}) also have the solution (\ref
{a2.6b}). If the vector $u^i$ is spacelike, $G_{ik}u^iu^k<0,$ then two
equations (\ref{a4.17}) have non-trivial solution, and the collinearity cone
does not degenerate into a line. The collinearity cone is a section of the
light cone $G_{ik}v^iv^k=0$ by the plane $G_{ik}u^iv^k=0$. If the vector $%
u^i $ is null, $G_{ik}u^iu^k=0,$ then equation (\ref{a4.18}) reduces to the
form
\begin{equation}
G_{ik}u^iu^k=0,\qquad G_{ik}u^iv^k=0  \label{a4.19}
\end{equation}
In this case (\ref{a2.6b}) is a solution, but besides there are spacelike
vectors $v^i$ which are orthogonal to null vector $u^i$ and the collinearity
cone does not degenerate into a line.

In the case of the proper $\sigma $-Riemannian space $G_{ik}u^iu^k>0,$ and
equation (\ref{a4.18}) reduces to the form
\begin{equation}
G_{ik}v^iv^k-\Delta _{ik}\left( \alpha u^i+v^i\right) \left( \alpha
u^k+v^k\right) =0  \label{a4.20}
\end{equation}
In this case $\Delta _{ik}\neq 0$ in general, and the collinearity cone does
not degenerate. $\Delta _{ik}$ depends on the curvature an on the distance
between the points $x$ and $x^{\prime }$. The more space curvature and the
distance $\rho (x,x^{\prime }),$ the more the collinearity cone aperture.

In the curved proper $\sigma $-Riemannian space there is an interesting
special case, when the collinearity cone degenerates . In any $\sigma $%
-Riemannian space the following equality takes place \cite{R62}
\begin{equation}
G_{ik}\sigma ^k=g_{ik}\sigma ^k,\qquad \sigma ^k\equiv g^{kl}\sigma _l
\label{a4.21}
\end{equation}
Then it follows from (\ref{a4.14}) that
\begin{equation}
\Delta _{ik}\sigma ^k=0  \label{a4.22}
\end{equation}
It means that in the case, when the vector $u^i$ is directed along the
geodesic, connecting points $x$ and $x^{\prime }$, i.e. $u^i=\beta \sigma ^i$%
, the equation (\ref{a4.20}) reduces to the form
\begin{equation}
\left( G_{ik}-\Delta _{ik}\right) v^iv^k=0,\qquad u^i=\beta \sigma ^i
\label{a4.23}
\end{equation}
If $\Delta _{ik}$ is small enough as compared with $G_{ik},$ then
eigenvalues of the matrix $G_{ik}-\Delta _{ik}$ have the same sign, as those
of the matrix $G_{ik}.$ In this case equation (\ref{a4.23}) has the only
solution (\ref{a2.6b}), and the collinearity cone degenerates.

\section{Discussion}

Thus, we see that in the $\sigma $-Riemannian geometry at the point $x$
there are many vectors parallel to given vector at the point $x^{\prime }$.
This set of parallel vectors is described by the collinearity cone.
Degeneration of the collinearity cone into a line, when there is only one
direction, parallel to the given direction, is an exception rather than a
rule, although in the proper Euclidean geometry this degeneration takes
place always. Nonuniformity of space destroys the collinearity cone
degeneration. In the proper Riemannian geometry, where the world function
satisfies the system (\ref{b2.39}), one succeeded in conserving this
degeneration for direction along the geodesic, connecting points $x$ and $%
x^{\prime }$. This circumstance is very important for degeneration of the
first order NGOs into geodesic, because degeneration of NGOs is connected
closely with the collinearity cone degeneration.

Indeed, definition of the first order tube (\ref{b1.3}), or (\ref{c3.2}) may
be written also in the form
\begin{equation}
\mathcal{T}\left( \mathcal{P}^1\right) \equiv \mathcal{T}_{P_0P_1}=\left\{
R\left| \; \mathbf{P}_0\mathbf{P}_1||\mathbf{P}_0\mathbf{R}\right. \right\}
,\qquad P_0,P_1,R\in \Omega ,  \label{b5.1}
\end{equation}
where collinearity $\mathbf{P}_0\mathbf{P}_1||\mathbf{P}_0\mathbf{R}$ of two
vectors $\mathbf{P}_0\mathbf{P}_1$ and $\mathbf{P}_0\mathbf{R}$ is defined
by the $\sigma $-immanent relation (\ref{a2.3}), which can be written in the
form

\begin{equation}
\mathbf{P}_0\mathbf{P}_1||\mathbf{P}_0\mathbf{R}:\quad F_2\left(
P_0,P_1,R\right) =\left|
\begin{array}{cc}
\left( \mathbf{P}_0\mathbf{P}_1.\mathbf{P}_0\mathbf{P}_1\right) & \left(
\mathbf{P}_0\mathbf{P}_1.\mathbf{P}_0\mathbf{R}\right) \\
\left( \mathbf{P}_0\mathbf{R}.\mathbf{P}_0\mathbf{P}_1\right) & \left(
\mathbf{P}_0\mathbf{R}.\mathbf{P}_0\mathbf{R}\right)
\end{array}
\right| =0  \label{b5.2}
\end{equation}
The form (\ref{b5.1}) of the first order tube definition allows one to
define the first order tube $\mathcal{T}(P_0,P_1;Q_0)$, passing through the
point $Q_0$ collinear to the given vector $\mathbf{P}_0\mathbf{P}_1$. This
definition has the $\sigma $-immanent form
\begin{equation}
\mathcal{T}(P_0,P_1;Q_0)=\left\{ R\left|\; \mathbf{P}_0\mathbf{P}_1||\mathbf{%
Q}_0\mathbf{R}\right. \right\} ,\qquad P_0,P_1,Q_0,R\in \Omega ,
\label{b5.3}
\end{equation}
where collinearity $\mathbf{P}_0\mathbf{P}_1||\mathbf{Q}_0\mathbf{R}$ of two
vectors $\mathbf{P}_0\mathbf{P}_1$ and $\mathbf{Q}_0\mathbf{R}$ is defined
by the $\sigma $-immanent relations (\ref{a2.3}), (\ref{a2.16b}). In the
proper Euclidean space the tube (\ref{b5.3}) degenerates into the straight
line, passing through the point $Q_0$ collinear to the given vector $\mathbf{%
P}_0\mathbf{P}_1$.

Let us define the set $\omega _{Q_0}=\left\{ \mathbf{Q}_0\mathbf{Q)|}Q\in
\Omega \right\} $ of vectors $\mathbf{Q}_0\mathbf{Q}$. Then
\begin{equation}
\mathcal{C}(P_0,P_1;Q_0)=\left\{ \mathbf{Q}_0\mathbf{Q|}Q\in \mathcal{T}%
(P_0,P_1;Q_0)\right\} \subset \omega _{Q_0}  \label{b5.4}
\end{equation}
is the collinearity cone of vectors $\mathbf{Q}_0\mathbf{Q}$ collinear to
vector $\mathbf{P}_0\mathbf{P}_1$. Thus, the one-dimensionality of the first
order tubes and the collinearity cone degeneration are connected phenomena.

In the Riemannian geometry the very special property of the proper Euclidean
geometry (the collinearity cone degeneration) is considered to be a property
of any geometry and extended to the case of Riemannian geometry. The line $%
\mathcal{L}$, defined as a continuous mapping (\ref{b5.5}) is considered to
be the most important geometrical object. This object is considered to be
more important, than the metric, and metric in the Riemannian geometry is
defined in terms of the shortest lines. Use of line as a basic concept of
geometry is inadequate for description of geometry and poses problems, which
appears to be artificial. For instance, the convexity problem, when
elimination of part of the point set $\Omega $ generates variation of
properties of other regions is a result of the metric definition via concept
of the line. Although choosing the world function in the proper way
(satisfying equations (\ref{b2.39})), one succeeded in conserving the
collinearity cone degeneration for geodesic lines, but for distant points $x$
and $x^{\prime }$ the collinearity cone does not degenerate, and the
absolute parallelism is absent in the Riemannian geometry. Instead of the
cone of collinear vectors one introduces concept of parallel transport of a
vector, where the result depends on the path of the transport. Practically,
it means that one vector of the vector cone is chosen and it is attributed
to some curve connecting the points $x$ and $x^{\prime }$.

Being a special case of T-geometry, the $\sigma $-Riemannian geometry does
not use the nonmetric concept of line at all. Here the nonmetric line is a
special geometrical object characteristic for the proper Euclidean geometry
which is a result of the collinearity cone degeneration. Instead of the
continuous mapping (\ref{b5.5}) one uses the mapping
\begin{equation}
m_n:\;\;I_n\rightarrow \Omega ,\qquad I_n=\left\{ 0,1,\ldots n\right\}
\subset \Bbb{Z}  \label{b5.6}
\end{equation}
which determines geometrical object $m_n$, called the $n$th order
multivector. \cite{R00}. The $n$th order multivector may be considered to be
some generalization of the $n$th order $\sigma $-subspace $M(\mathcal{P}^n)$%
, and definition (\ref{b5.6}) of multivector appears to be $\sigma $%
-immanent. Application of mappings (\ref{b5.6}) is sufficient for
description of any geometry, because all geometric objects are determined as
subsets of the space $\Omega $ (not as mappings). Use of such complicated
mappings as (\ref{b5.5}) is not necessary. For instance, to investigate the
properties of the first order tube $\mathcal{T}_{P_0P_1}\subset \Omega $
(geodesic), one needs to investigate the set $\mathcal{T}=\left\{
P_0\right\} \otimes \left\{ P_1\right\} \otimes \mathcal{T}_{P_0P_1}\subset
\Omega ^3,$ satisfying the condition $F_3(\mathcal{T})=0.$ Here the mapping $%
F_3$ is known and fixed. Only zeros of the function $F_3$, having the form $%
\mathcal{T}=\left\{ P_0\right\} \otimes \left\{ P_1\right\} \otimes \mathcal{%
T}_{P_0P_1}$, are investigated. Power of the set $\mathcal{T}$ is much less
than the power of the set of all mappings (\ref{b5.5}), and investigation of
$\mathcal{T}$ is not so complicated as investigation of mappings (\ref{b5.5}%
).

One can reduce the power of the set of all mappings (\ref{b5.5}), imposing
some additional restrictions on mapping (\ref{b5.5}), but nothing can change
the fact that the mapping (\ref{b5.5}) is an attribute of the proper
Euclidean geometry and is not an attribute of a geometry in itself. The
convexity problem confirms this. The real space-time may appear not to have
property of the collinearity cone degeneration \cite{R91}. Insisting on the
mapping (\ref{b5.5}) as the main tool of geometry investigation, one closes
the door for real investigations of geometry and shows a wrong way for them.

Besides purely logical arguments in favour of the T-geometry approach there
are arguments of applied character. The fact is that application of
T-geometry to the space-time model construction leads to new encouraging
results \cite{R90,R91}. Consideration of uniform isotropic continuous model
with zeroth curvature leads to a class of models, distinguishing by the
shape of the tube. This class contains the well known Minkowski model, for
which the timelike tubes degenerate into lines and which is not optimal,
because it does not enable to describe quantum phenomena without using the
quantum principles. Other (nondegenerate) models of this class have the
following properties: (1) geometrization of mass of a particle described by
the broken tube (\ref{a1.20}), (2) stochasticity of the world tube of a free
particle which is conditioned by the collinearity cone non-degeneracy.

It turns out that it is possible one to choose optimal space-time model, for
which the statistical description of stochastic free particle tubes
coincides with the quantum description in terms of the Schr\"odinger
equation. The quantum constant $\hbar$ appears to be a space-time property,
introducing some "elementary length" (it is connected with the thickness of
the particle world tube). As a result one does not need the quantum
principles, and the quantum theory looks as a conception, created for
compensation of our incorrect ideas on the space-time geometry at small
distances.

\end{document}